\documentclass[12pt]{amsart}
\usepackage{amsmath,amssymb,amsthm,upref,graphicx,mathrsfs}

\usepackage{color}
\usepackage[
  colorlinks=true,
  linkcolor=blue,
  citecolor=blue,
  urlcolor=blue]{hyperref}

\numberwithin{equation}{section}

\newtheorem{theorem}{Theorem}[section]

\newtheorem{lemma}[theorem]{Lemma}

\theoremstyle{definition}

\newtheorem{remark}[theorem]{Remark}

\newcommand{\Mm}{\mathcal{M}}

\newcommand{\Rn}{\mathbb{R}^n}
\newcommand{\Z}{\mathbb{Z}}

\let \a=\alpha

\let \O=\Omega

\begin{document}

\title[Weighted estimates for the multisublinear maximal function]
  {Weighted estimates for the multisublinear maximal function}

\authors

\author[W. Chen]{Wei Chen}
\address{Wei Chen \\ School of Mathematical Sciences,
Yangzhou University, 225002 Yangzhou, China}
\email{weichen@yzu.edu.cn}

\author[W. Dami\'an]{Wendol\'in Dami\'an}
\address{Wendol\'in Dami\'an\\
Departamento de An\'alisis Matem\'atico, Facultad de Matem\'aticas,
Universidad de Sevilla, 41080 Sevilla, Spain}
\email{wdamian@us.es}

\makeatletter
\renewcommand{\@makefntext}[1]{#1}
\makeatother \footnotetext{\noindent
 The first author is supported by the National Natural Science Foundation of China (Grant No. 11101353), the Natural Science Foundation of Jiangsu Education Committee (Grant No. 11KJB110018) and the Natural Science Foundation of Jiangsu Province (Grant No. BK2012682). The second author is supported by Junta de Andaluc\'ia (Grant No. P09-FQM-4745)}


%
%
\begin{abstract}
 A formulation of the Carleson embedding theorem in the multilinear setting is proved which allows  to obtain a multilinear analogue of Sawyer's two weight theorem for the multisublinear maximal function $\Mm$ introduced in \cite{LOPTT}.  A multilinear version of the $B_p$ theorem from \cite{HP}
 is also obtained and a mixed $A_{\vec P}-W_{\vec P}^{\infty}$ bound for $\Mm$
 is proved as well.
\end{abstract}

\maketitle

%
%
\section{Introduction}

The beginning of the modern theory of weights was originated in the works of R. Hunt, B. Muckenhoupt, R. Wheeden, R. Coifman and C. Fefferman in the decade of the 70's. In \cite{Mu} B. Muckenhoupt characterized the class of weights $u,v$ for which the following weak inequality holds

\begin{equation}\label{weak_two_weight_problem}
  \sup_{\lambda>0} \lambda^p \int_{\{Mf > \lambda\}} u(x)dx \leq C \int_{\Rn} |f(x)|^p v(x) dx, \, f\in L^p(v),
\end{equation}
where $M$ denotes the Hardy--Littlewood maximal operator and $p\geq 1$. This condition on the weights is known as $A_p$ condition, namely

\begin{equation*}
  [u,v]_{A_p} := \sup_{Q} \left(\frac{1}{|Q|}\int_Q u(x) dx\right) \left(\frac{1}{|Q|}\int_Q v(x)^{-\frac{1}{p-1}}\right)^{p-1} < \infty, \, p>1
\end{equation*}
where the supremum is taken over all the cubes in $\Rn$. When $p=1$, the term $(\int_Q \frac{1}{|Q|}v(x)^{-\frac{1}{p-1}})^{p-1}$
must be understood as $(\inf_Q v)^{-1}$.  In the particular case when $u=v$, Muckenhoupt also proved that the following strong estimate

\begin{equation*}\label{strong_one_weight_problem}
  \int_{\Rn} (Mf(x))^p v(x) dx \leq C \int_{\Rn} |f(x)|^p v(x) dx, \, f\in L^p(v),
\end{equation*}
holds if and only if $v$ satisfies the $A_p$ condition. However,  the problem of finding a condition on the weights $u,v$ satisfying the strong estimate above was more complicated. In \cite{Saw} E. Sawyer characterized the two weight inequality, showing that $M:L^p(v)\longrightarrow L^p(u)$ if and only if the pair $(u,v)$ satisfies the following testing condition known as Sawyer's $S_p$ condition

\begin{equation}\label{Sp_constant_linear}
  [u,v]_{S_p} = \sup_{Q}  \left(\frac{\int_{Q} M(\chi_Q \sigma)^p u dx}{\sigma(Q)}\right)^{1/p} < \infty,
\end{equation}
where $\sigma=v^{1-p'}$ and $1<p<\infty$. Motivated by these results the theory of weighted inequalities developed rapidly, not only for the Hardy--Littlewood maximal operator but also for some of the main operators in Harmonic Analysis like Calder\'on--Zygmund operators.  Much later the interest focused in determining the sharp dependence of the $L^p(w)$ operator norm in term of the relevant  constant involving the weights.

\par On this point, the problem for the Hardy--Littlewood maximal operator was solved by S. Buckley \cite{Bu} who proved

\begin{equation}\label{Buckley_thm}
  ||M||_{L^p(w)} \leq C \,p'\,[w]_{A_p}^{\frac{1}{p-1}},
\end{equation}
where $C$ is a dimensional constant.  Motivated by this result and others, K. Moen found in \cite{Moen_S} a quantitative form of E. Sawyer's result mentioned above in terms of Sawyer's condition \eqref{Sp_constant_linear},  namely
\begin{equation}\label{Sp_sharp_result_Moen}
  ||M||_{L^p(v)\longrightarrow L^p(u)} \approx [u,v]_{S_p}.
\end{equation}

Recently,  T. Hyt\"onen and C. P\'erez in \cite{HP} (see also \cite{HPR} for a better result and a simplified proof) improved Buckley's bound \eqref{Buckley_thm} replacing a portion of the $A_p$ constant by the weaker $A_{\infty}$ constant as defined by Fujii in \cite{Fu} and later used in the work of J.M. Wilson \cite{Wil}. The $A_{\infty}$ constant is defined as follows

\begin{equation}\label{Fujii_Wilson_linear_constant}
  [w]_{A_{\infty}} := \sup_Q \frac{1}{w(Q)}\int_{Q} M(w\chi_Q),
\end{equation}
where the supremum is taken over all the cubes $Q$ in $\Rn$. In \cite{HP} the authors show in a two-weight setting and for $p>1$ that

\begin{equation}\label{HP_Bp_thm}
  ||M(f\sigma)||_{L^p(w)} \leq C p' (B_p[w,\sigma])^{1/p} ||f||_{L^p(\sigma)},
\end{equation}
and

\begin{equation}\label{HP_Ap_Ainfty_thm}
  ||M(f\sigma)||_{L^p(w)} \leq C p' ([w]_{A_p} [\sigma]_{A_{\infty}})^{1/p} ||f||_{L^p(\sigma)},
\end{equation}
where $C$ in both inequalities is a dimensional constant and

\begin{equation}\label{Bp_constant_linear}
  B_p[w,\sigma]:= \sup_Q \left(\frac{1}{|Q|}\int_Q w\right) \left(\frac{1}{|Q|} \int_Q \sigma \right)^p \exp{\left(\frac{1}{|Q|}\int_Q \log \sigma^{-1} \right)},
\end{equation}
is known as the $B_p$ constant of the weights $w$ and $\sigma$. This constant clearly satisfies

\begin{equation*}
  [w]_{A_p} \leq B_p[w,\sigma] \leq [w]_{A_p} [\sigma]_{A_{\infty}}',
\end{equation*}
where $[\sigma]_{A_{\infty}}'$ denotes the $A_{\infty}$ constant introduced by S. Hrus\v{c}\v{e}v  in \cite{Hru} defined as follows

\begin{equation*}
  [w]_{A_{\infty}}' = \sup_{Q} \left(\frac{1}{|Q|}\int_{Q} w \right) \exp{\left(\frac{1}{|Q|}\int_{Q}\log w^{-1}\right)}.
\end{equation*}
In the one weight setting and by a standard change-of-weight argument, \eqref{HP_Ap_Ainfty_thm} implies

\begin{equation}\label{HP_Ap_Ainfty_thm_one_weight}
  ||M||_{L^p(w)} \leq C p' ([w]_{A_p} [\sigma]_{A_{\infty}})^{1/p},
\end{equation}
where $\sigma=w^{1-p'}$.

\par The aim of this article is to give some multilinear analogues of some of the above mentioned results following the spirit of the theory of multiple weights developed in \cite{LOPTT}.

\par The paper is organized as follows. Some preliminary definitions and results are summarized in Sect.~\ref{sec:2}. In Sect.~\ref{sec:3} we give the statements of the main results on this paper and some remarks on them. Finally, in Sect.~\ref{sec:4} we give all the proofs of our results.

Throughout this paper, we will use the notation $A\lesssim B$ to indicate that there is a constant $c$, independent of the weight constant, such that $A\leq c B$.

%
%

\section{Preliminaries}
\label{sec:2}

Before stating and proving our main results, we first recall some basics related to the theory of multilinear weighted inequalities as well as we introduce the definition of some constants involved in the multiple theory of weights and dyadic grids.

\subsection{Some basics on multilinear weighted inequalities}
\label{subsec:2_1}

One of the main objects of the theory of multiple weights is the following extension of the classical Hardy--Littlewood maximal function. Given $\overrightarrow f=(f_1,\dots,f_m)$, we define following \cite{LOPTT} the multi(sub)linear maximal operator $\mathcal M$ by
$$\mathcal M(\overrightarrow f\,)(x)=\sup_{Q\ni x}\prod_{i=1}^m\frac{1}{|Q|}\int_Q|f_i(y_i)|dy_i,$$
where the supremum is taken over all cubes $Q$ containing $x$. The importance of this operator stems from the fact that it controls the class of multilinear Calder\'on--Zygmund operators as it is shown in \cite{LOPTT}. A particular example of this relationship is the class of weights characterizing the weighted $L^p$ spaces for which both operators are bounded. To define this class of weights we let  $\overrightarrow w=(w_1,\dots,w_m)$ and $\overrightarrow P=(p_1,\dots, p_m)$ such that $1<p_1,\ldots,p_m<\infty$. Set $\frac{1}{p}=\frac{1}{p_1}+\dots+\frac{1}{p_m}$ and $\nu_{\overrightarrow w}=\prod_{i=1}^mw_i^{p/p_i}$. We say that $\overrightarrow w$ satisfies the $A_{\overrightarrow P}$ condition if

$$[\overrightarrow w]_{A_{\overrightarrow P}}=\sup_{Q}\Big(\frac{1}{|Q|}\int_Q\nu_{\overrightarrow
w}\Big)\prod_{i=1}^m\Big(\frac{1}{|Q|}\int_Q
w_i^{1-p'_i}\Big)^{p/p'_i}<\infty.
$$

It is easy to see that in the linear case (that is, if $m=1$) $[\overrightarrow w]_{A_{\overrightarrow P}}=[w]_{A_p}$ is the usual $A_p$ constant.

In \cite{LOPTT} the following multilinear extension of the Muckenhoupt $A_p$ theorem for the maximal function was obtained:  the inequality

\begin{equation}\label{strong}
\|{\mathcal M}(\overrightarrow f\,)\|_{L^{p}(\nu_{\overrightarrow w})}\le C \prod_{i=1}^m\|f_i\|_{L^{p_i}(w_i)}
\end{equation}
holds for every  $\overrightarrow f$ if and only if $\overrightarrow w$ satisfies the $A_{\overrightarrow P}$ condition.

\par Very recently in \cite{DLP} A. Lerner, C. P\'erez and the second author proved a multilinear version of Buckley's result as well as a full analogue of \eqref{HP_Ap_Ainfty_thm}. In this work the authors found that the multilinear version of \eqref{HP_Ap_Ainfty_thm} is sharp when $m\geq 1$, although is much more complicated to do the same for Buckley's result. In this case, several partial results were obtained  in \cite{DLP}  which have been improved in \cite{LMS}.

\subsection{Some constants on multiple weight theory.}
\label{subsec:2_2}

Next we state the notation that we will follow in the sequel related to some constants involved in the multiple theory of weights. To define these constants, let $w_1,\ldots,w_m$ and $v$ be weights and let us denote $\overrightarrow{w}=(w_1,\ldots,w_m)$.  Also let $1<p_1,\ldots,p_m<\infty$ and $p$ be numbers such that $\frac{1}{p}=\frac{1}{p_1}+\dots+\frac{1}{p_m}$ and denote $\overrightarrow P = (p_1,\ldots, p_m)$.

\par We say that $(v,\overrightarrow w)$ satisfies the $A_{\overrightarrow P}$ condition if
  \begin{equation}\label{Ap_constant}
  [v,\overrightarrow w]_{A_{\overrightarrow P}}:=\sup_{Q}\Big(\frac{1}{|Q|}\int_Q v \Big)\prod_{i=1}^m\Big(\frac{1}{|Q|}\int_Q w_i^{1-p'_i}\Big)^{p/p'_i}<\infty.
  \end{equation}

  In particular when $v=\nu_{\overrightarrow w}:=\prod_{i=1}^m w_i^{\frac{p}{p_i}}$, we will write $[\nu_{\overrightarrow w},\overrightarrow w]_{A_{\overrightarrow P}}$ as $[\overrightarrow w]_{A_{\overrightarrow P}}$.

\par Next we define the multilinear analogues of the $A_\infty$ constant defined by Fujii in \cite{Fu}, the $B_p$ constant defined by Hyt\"onen and P\'erez in \cite{HP} and the $S_p$ constant defined by Sawyer in \cite{Saw}, respectively. We say that

\begin{enumerate}
\item $\overrightarrow w$ satisfies the $W_{\overrightarrow P}^\infty$ condition if
    \begin{equation*}\label{Fujii_constant}
      [\overrightarrow w]_{W_{\overrightarrow P}^\infty}=\sup_Q \Big(\int_Q\prod^m_{i=1}M(w_i\chi_Q)^{\frac{p}{p_i}}dx\Big)\Big(\int_Q\prod^m_{i=1}w_i^{\frac{p}{p_i}} dx\Big)^{-1}<\infty.
    \end{equation*}

\item $(v,\overrightarrow w)$ satisfies the $B_{\overrightarrow{P}}$ condition if
    \begin{equation*}\label{Bp_constant}
    [v,\overrightarrow w]_{B_{\overrightarrow P}}:= \sup_Q \frac{v(Q)}{|Q|} \Big(\prod^m_{i=1}\frac{w_i(Q)}{|Q|}\Big)^p
  \exp\Big(\frac{1}{|Q|}\int_Q\log\prod^m_{i=1}w_i^{-\frac{p}{p_i}} dx\Big)<\infty.
    \end{equation*}

\item $(v,\overrightarrow w)$ satisfies the $S_{\overrightarrow{P}}$ condition if
   \begin{equation*}\label{Sp_constant}
     [v,\overrightarrow w]_{S_{\overrightarrow P}} = \sup_Q \Big(\int_Q \Mm(\overrightarrow{\sigma\chi_Q})^{p} v dx \Big)^{\frac{1}{p}}  \Big(\prod^m_{i=1}\sigma_i(Q)^{\frac{1}{p_i}}\Big)^{-1}<\infty,
   \end{equation*}
where $\overrightarrow{\sigma\chi_Q}=(\sigma_1\chi_Q,\ldots,\sigma_m\chi_Q)$ and $\sigma_i=w_i^{1-p_i'}$ for all $i=1,\ldots,m$ and all the suprema in the above definitions are taken over all cubes $Q$ in $\Rn$.
\end{enumerate}

\noindent Additionally we define a multiple Reverse H\"older condition that we will use in the following. We say that $\overrightarrow w$ satisfies the $RH_{\overrightarrow P}$ condition if
there exists a positive constant $C$ such that

  \begin{equation}\label{RH_constant}
     \prod^m_{i=1}\Big(\int_Q\sigma_i dx \Big)^{\frac{p}{p_i}} \leq  C \int_Q\prod^m_{i=1}\sigma_i^{\frac{p}{p_i}} dx,
  \end{equation}
where $\sigma_i=w_i^{1-p_i'}$ for $i=1,\ldots,m$. We denote by $[\overrightarrow w]_{RH_{\overrightarrow P}}$ the smallest constant $C$ in \eqref{RH_constant}.

\subsection{Dyadic grids.}
\label{subsec:2_3}

Recall that the standard dyadic grid ${\mathcal D}$ in ${\mathbb R}^n$ consists of the cubes

\begin{equation*}
  2^{-k}([0,1)^n+j),\quad k\in{\mathbb Z}, j\in{\mathbb Z}^n.
\end{equation*}

By a {\it general dyadic grid} ${\mathscr{D}}$ we mean a collection of cubes with the following properties:

\begin{enumerate}
  \item For any $Q\in {\mathscr{D}}$ its sidelength $\ell_Q$ is $2^k, k\in {\mathbb Z}$
  \item $Q\cap R\in\{Q,R,\emptyset\}$ for any $Q,R\in {\mathscr{D}}$.
  \item The cubes of a fixed sidelength $2^k$ form a partition of ${\mathbb R}^n$.
\end{enumerate}
We say that $\{Q_j^k\}$ is a {\it sparse family} of cubes if:

\begin{enumerate}
  \item The cubes $Q_j^k$ are disjoint in $j$, with $k$ fixed.
  \item If $\Omega_k=\cup_jQ_j^k$, then $\Omega_{k+1}\subset~\Omega_k$.
  \item $|\Omega_{k+1}\cap Q_j^k|\le \frac{1}{2}|Q_j^k|$.
\end{enumerate}

With each sparse family $\{Q_j^k\}$ we associate the sets $E_j^k=Q_j^k\setminus \O_{k+1}$.  Observe that the sets $E_j^k$ are pairwise disjoint and $|Q_j^k|\le 2|E_j^k|$.

In the sequel we will use the following lemmas that could be found in \cite{HP} and \cite{DLP}, respectively.

\begin{lemma}\label{prhp} There are $2^n$ dyadic grids ${\mathscr{D}}_{\a}$ such that for any cube $Q\subset {\mathbb R}^n$ there exists a cube $Q_{\a}\in {\mathscr{D}}_{\a}$
such that $Q\subset Q_{\a}$ and $\ell_{Q_{\a}}\le 6\ell_Q$.
\end{lemma}

\begin{lemma}\label{appr}
For any non-negative integrable $f_i,i=1,\dots,m$, there exist sparse families ${\mathcal S}_{\a}\in {\mathscr{D}}_{\a}$ such that for all $x\in {\mathbb R}^n$,
$${\mathcal M}(\overrightarrow f\,)(x)\le (2\cdot 12^n)^{m}\sum_{\a=1}^{2^n}
{\mathcal A}_{{\mathscr{D}_{\a}},{\mathcal S}_{\a}}(\overrightarrow f\,)(x),$$
where $\overrightarrow{f}=(f_1,\ldots,f_m)$ and given a sparse family $\mathcal S =\{Q_j^k\}$ of cubes from a dyadic grid $\mathscr{D}$, the operator ${\mathcal A}_{{\mathscr{D}},{\mathcal S}}$ is given by

$${\mathcal A}_{{\mathscr{D}},{\mathcal S}}(\overrightarrow f\,) = \sum_{j,k} \left(\prod_{i=1}^m (f_i)_{Q_j^k}\right)\chi_{Q_j^k}.$$

\end{lemma}

\section{Main results}
\label{sec:3}

\par In this section we summarize the main results on this work. Firstly we state the main tool of this paper. This lemma  extends to the multilinear setting a nonstandard formulation of the (dyadic) Carleson embedding theorem proved in \cite{HP} and it will allow us to prove our main results.

\begin{lemma}\label{Carleson_lemma_multi} Suppose that the nonnegative numbers $\{a_Q\}_Q$ satisfy

\begin{equation}\label{Carleson_assumption}
 \sum_{Q\subset R} a_Q \leq A \int_R \prod_{i=1}^m \sigma_i^{\frac{p}{p_i}}dx, \, \forall R \in \mathscr{D}
\end{equation}

\noindent where $\sigma_i$ are weights for $i=1,\ldots,m$. Then for all $1 < p_i < \infty$  and $p\in (1,\infty)$ satisfying $\frac{1}{p}=\frac{1}{p_1}+\dots+\frac{1}{p_m}$ and for all $f_i\in L^{p_i}(\sigma_i)$,

\begin{equation}\label{Carleson}\begin{split}
\left(\sum_{Q\in \mathscr{D}} a_Q \Big(\prod_{i=1}^{m} \frac{1}{\sigma_i(Q)}\int_{Q}f_i(y_i)\sigma_i(y_i)d y_i\Big)^p\right)^{1/p} \leq& A ||\Mm^d_{\overrightarrow{\sigma}}(\overrightarrow f)||_{L^p(\nu_{\overrightarrow \sigma})} \\
\leq& A \prod_{i=1}^{m} p'_i ||f_i||_{L^{p_i}(\sigma_i)},
\end{split}\end{equation}
where $\Mm^d_{\overrightarrow{\sigma}}(\overrightarrow{f})=\displaystyle{\sup_{\substack{Q\ni x \\ Q\in\mathscr{D}}}}\prod_{i=1}^m \frac{1}{\sigma_i(Q)}\int_Q |f_i(y_i)|\sigma_i(y_i)dy_i$.

\end{lemma}

Next we establish a generalization of Sawyer's theorem to the multilinear setting. Very recently it was shown in \cite{LXY} a multilinear version of Sawyer's theorem using a kind of monotone property on the weights. We establish here another condition that is a sort of reverse H\"older inequality in the multilinear setting (see Sect.~\ref{sec:2} for definition) and that was used by the first author in \cite{CL} in the setting of martingale spaces. When $m=1$ this reverse H\"older condition is superfluous and we recover the linear result of Moen \eqref{Sp_sharp_result_Moen}.

\begin{theorem}\label{Sp_theorem}
Let $1< p_i < \infty$, $i=1,\ldots,m$ and $\frac{1}{p}=\frac{1}{p_1}+\ldots+\frac{1}{p_m}$.  Let $v$ and $w_i$ be weights. If we suppose that $\overrightarrow w \in RH_{\overrightarrow P}$ then there exists a positive constant $C$ such that

    \begin{equation}\label{eqn_Sp}
      ||\Mm(\overrightarrow{f\sigma})||_{L^p(v)} \leq C \prod_{i=1}^m ||f_i||_{L^{p_i}(\sigma_i)}, ~f_i\in L^{p_i}(\sigma_i),
    \end{equation}
where $\sigma_i=w_i^{1-p_i'}$, if and only if $(v,\overrightarrow w)\in S_{\overrightarrow P}$.  Moreover, if we denote the smallest constant $C$ in \eqref{eqn_Sp} by $||\Mm||$, we obtain

\begin{equation}\label{relationship_Sp}
  [v,\overrightarrow w]_{S_{\overrightarrow P}} \lesssim ||\Mm|| \lesssim [v,\overrightarrow w]_{S_{\overrightarrow P}} [\overrightarrow w]^{1/p}_{RH_{\overrightarrow P}}.
\end{equation}

\end{theorem}

\noindent Here we make some remarks related to the previous theorem.

\begin{remark}
  In the particular case when $v=\nu_{\overrightarrow w}$, the following statements are equivalent:

  \begin{enumerate}
    \item $\overrightarrow w \in A_{\overrightarrow P}$.
    \item $\sigma_i=w_i^{1-p_i'}\in A_{mp_i'}$, for $i=1,\ldots,m$ and $\nu_{\overrightarrow w}\in A_{mp}$.
    \item $(\nu_{\overrightarrow w},\overrightarrow w)\in S_{\overrightarrow P}$.
    \item There exists a positive constant $C$ such that
      \begin{equation}\label{eqn_rel_S_p}
        ||\Mm(\overrightarrow f)||_{L^p(\nu_{\overrightarrow w})} \leq C \prod_{i=1}^m ||f_i||_{L^{p_i}(w_i)}, ~f_i \in L^{p_i}(w_i).
      \end{equation}
  \end{enumerate}

  Indeed, the equivalence between $1.$, $2.$ and $4.$ was proved in \cite[Th. 3.6, Th. 3.7]{LOPTT}. It can be easily seen that in this particular case $[\nu_{\overrightarrow w},\overrightarrow w]_{S_{\overrightarrow P}}\lesssim ||\Mm||$ where $||\Mm||$ denotes the smallest constant in \eqref{eqn_rel_S_p} and $[\overrightarrow w]_{A_{\overrightarrow P}}\lesssim [\nu_{\overrightarrow w},\overrightarrow w]_{S_{\overrightarrow P}}^p$. Therefore we have that $4.$ implies $3.$ and $3.$ implies $1.$. So we have obtained that all the statements are equivalent.

  \par Additionally, following \cite[Th. 1.1]{DLP}, we also have that $||\Mm|| \lesssim [\overrightarrow w]_{A_{\overrightarrow P}}^{1/p}\prod_{i=1}^m [\sigma_i]_{\infty}^{\frac{1}{p_i}}$. So, we have obtained

  \begin{equation}\label{relationship_Sp}
  [\overrightarrow w]_{A_{\overrightarrow P}}^{1/p} \lesssim [v_{\overrightarrow w},\overrightarrow w]_{S_{\overrightarrow P}} \lesssim ||\Mm|| \lesssim [\overrightarrow w]_{A_{\overrightarrow P}}^{1/p}\prod_{i=1}^m [\sigma_i]_{\infty}^{\frac{1}{p_i}}.
  \end{equation}

  \end{remark}

  \begin{remark}
    As we have observed in the previous remark, $RH_{\overrightarrow P}$ condition is not necessary when $v=\nu_{\overrightarrow w}$ in Theorem \ref{Sp_theorem}. We are not sure if this condition can be removed in the general case.
  \end{remark}

Making use of the analogue of the $B_p$ constant within the multilinear setting already defined in Sect.~\ref{sec:2}, we obtain an extension of \eqref{HP_Bp_thm}.

\begin{theorem}\label{Bp_theorem}
  Let $1<p_i<\infty$, $i=1,\ldots,m$ and $\frac{1}{p}=\frac{1}{p_1}+\ldots+\frac{1}{p_m}$.  Let $v$ and $w_i$ be weights. Then

  \begin{equation}\label{Bp_eqn}
     ||\Mm(\overrightarrow{f\sigma})||_{L^p(v)} \lesssim [v,\overrightarrow \sigma]_{B_{\overrightarrow P}}^{1/p} \prod_{i=1}^m ||f_i||_{L^{p_i}(\sigma_i)}, ~f_i\in L^{p_i}(\sigma_i),
  \end{equation}
  where $\sigma_i=w_i^{1-p_i'}$, $\overrightarrow{\sigma}=(\sigma_1,\ldots,\sigma_m)$ and $\overrightarrow{f\sigma}=(f_1\sigma_1,\ldots,f_m\sigma_m)$.
\end{theorem}

And finally, using the generalization of the Fujii--Wilson $A_{\infty}$ constant $[\overrightarrow w]_{W_{\overrightarrow P}^{\infty}}$ and the two-weight constant $[v,\overrightarrow w]_{A_{\overrightarrow P}}$ defined in Sect.~\ref{sec:2}, we get a mixed $A_{\overrightarrow P}-W_{\overrightarrow P}^{\infty}$ bound for $\Mm$ that extends \eqref{HP_Ap_Ainfty_thm} to the multilinear setting.

\begin{theorem}\label{Mixed_Ap_Ainfty_theorem}

    Let $1<p_i<\infty$, $i=1,\ldots,m$ and $\frac{1}{p}=\frac{1}{p_1}+\ldots+\frac{1}{p_m}$.  Let $v$ and $w_i$ be weights. Then

    \begin{equation}\label{Mixed_eqn}
      ||\Mm(\overrightarrow{f\sigma})||_{L^p(v)} \lesssim ([v,\overrightarrow w]_{A_{\overrightarrow P}}[\overrightarrow{\sigma}]_{W^{\infty}_{\overrightarrow P}})^{1/p} \prod_{i=1}^m ||f_i||_{L^{p_i}(\sigma_i)}, ~f_i\in L^{p_i}(\sigma_i),
    \end{equation}

    where $\sigma_i=w_i^{1-p_i'}$, $\overrightarrow{\sigma}=(\sigma_1,\ldots,\sigma_m)$ and $\overrightarrow{f\sigma}=(f_1\sigma_1,\ldots,f_m\sigma_m)$.
\end{theorem}

\section{Proofs}
\label{sec:4}


We start proving Lemma \ref{Carleson_lemma_multi}, that is, the multilinear version of the dyadic Carleson embedding theorem. It follows a scheme of proof similar to the one used by Hyt\"onen and P\'erez in \cite{HP}.

\begin{proof}[Lemma \ref{Carleson_lemma_multi}]
  Let us see the sum
    $$\sum_{Q\in\mathscr{D}} a_Q \left(\prod_{i=1}^m \frac{1}{\sigma_i(Q)}\int_Q f_i(y_i)\sigma_i(y_i)dy_i \right)^p$$
  as an integral on a measure space $(\mathscr{D},2^{\mathscr{D}},\mu)$ built over the set of dyadic cubes $\mathscr{D}$, assigning to each $Q\in\mathscr{D}$ the measure $a_Q$. Thus

  \begin{equation*}\begin{split}
  &\sum_{Q\in\mathscr{D}} a_Q \left(\prod_{i=1}^m \frac{1}{\sigma_i(Q)}\int_Q f_i(y_i)\sigma_i(y_i)dy_i \right)^p = \\ &=
  \int_0^{\infty} p \lambda^{p-1} \mu\left\{Q\in\mathscr{D}: \prod_{i=1}^m \frac{1}{\sigma_i(Q)} \int_Q f_i(y_i)\sigma_i(y_i)dy_i>\lambda\right\} \\ &=: \int_{0}^{\infty} p \lambda^{p-1}\mu(\mathscr{D}_{\lambda})d\lambda.
  \end{split}\end{equation*}

  Let us denote by $\mathscr{D}_{\lambda}^{*}$ the set of maximal dyadic cubes $R$ with the property that $\prod_{i=1}^m \frac{1}{\sigma_i(Q)}\int_R f_i(y_i)\sigma_i(y_i)dy_i >\lambda$. Then the cubes $R\in\mathscr{D}_{\lambda}^*$ are disjoint and their union is equal to the set $\{\Mm^d_{\overrightarrow{\sigma}}(\overrightarrow f)>\lambda\}$. Thus

  \begin{equation*}\begin{split}
    \mu(\mathscr{D}_{\lambda}) &= \sum_{Q\in\mathscr{D}_{\lambda}} a_Q \leq \sum_{R\in\mathscr{D}_{\lambda}^{*}}\sum_{Q\subset R} a_Q \\
    &\leq A \sum_{R\in\mathscr{D}^{*}_{\lambda}} \int_R \prod_{i=1}^m \sigma_i^{\frac{p}{p_i}} dx \\
    &= A \int_{\{\Mm_{\overrightarrow{\sigma}}^d(\overrightarrow f)>\lambda\}} \prod_{i=1}^m \sigma_i^{\frac{p}{p_i}}dx.
  \end{split}\end{equation*}
  Then we obtain

  \begin{equation*}\begin{split}
    \sum_{Q\in\mathscr{D}} a_Q \left(\prod_{i=1}^m \frac{1}{\sigma_i(Q)}\int_Q f_i(y_i)\sigma_i(y_i)dy_i \right)^p &\leq A\int_{0}^\infty p \lambda^{p-1} \int_{\{\Mm_{\overrightarrow\sigma}^d(\overrightarrow f)>\lambda\}} \prod_{i=1}^m \sigma_i^{\frac{p}{p_i}} dx d\lambda \\
    &= A \int_{\Rn} \Mm_{\overrightarrow{\sigma}}^d(\overrightarrow f)^p \prod_{i=1}^m \sigma_i^{\frac{p}{p_i}}dx \\
    &\leq A \int_{\Rn} \prod_{i=1}^m ((M_{\sigma_i}^d(f_i))^{p_i}\sigma_i)^{\frac{p}{p_i}}dx \\
    &\leq A \prod_{i=1}^m \left(\int_{\Rn} (M_{\sigma_i}^d(f_i))^{p_i} \sigma_i dx \right)^{\frac{p}{p_i}} \\
    &\leq A \prod_{i=1}^m \left(p_i'\right)^p \left(\int_{\Rn}|f_i|^{p_i}\sigma_i dx \right)^{\frac{p}{p_i}},
  \end{split}\end{equation*}
  where we have used that $\Mm_{\overrightarrow{\sigma}}^d(\overrightarrow f)\leq\prod_{i=1}^m M_{\sigma_i}^d(f_i)$, H\"older's inequality and the boundedness properties of $M_{\sigma_i}^d(f_i)$ in $L^{p_i}(\sigma_i)$.
\end{proof}


Next we prove Theorem \ref{Sp_theorem} making use of Lemma \ref{Carleson_lemma_multi}.

\begin{proof}[Theorem \ref{Sp_theorem}]
  It is clear that \eqref{eqn_Sp} implies the $S_{\overrightarrow P}$ condition without using that $(v,\overrightarrow w)\in RH_{\overrightarrow P}$. Thus, it remains to prove that $(v,\overrightarrow w)\in S_{\overrightarrow P}$ implies \eqref{eqn_Sp} to complete the proof of the theorem.

  \par By Lemma \ref{prhp}, it suffices to prove the theorem for the dyadic maximal operators ${\mathcal M}^{{\mathscr{D}}_{\a}}$. Since the proof is independent of the particular dyadic grid, without loss of generality we consider ${\mathcal M}^{d}$ taken with respect to the standard dyadic grid ${\mathcal D}$. Next we proceed as in the proof of Lemma \ref{appr}.  Let $a=2^{m(n+1)}$ and for $k\in\Z$ consider the following sets

  \begin{equation*}
    \O_k = \{ x\in\Rn: \Mm^d(\overrightarrow{f\sigma}) > a^k\}.
  \end{equation*}

  Then we have that $\O_k=\cup_j Q_{j}^{k}$, where the cubes $Q_j^k$ are pairwise disjoint with $k$ fixed, and

  \begin{equation*}
    a^k < \prod_{i=1}^m \frac{1}{|Q_j^k|}\int_{Q_j^k} |f_i(y_i)|\sigma_i(y_i)dy_i \leq 2^{mn}a^k.
  \end{equation*}

  It follows that

  \begin{equation*}\begin{split}
    \int_{\Rn} \Mm^d(\overrightarrow{f\sigma})^p v dx &= \sum_{k} \int_{\O_k\setminus\O_{k+1}} \Mm^d(\overrightarrow{f\sigma})^p v dx \\
    &\leq a^p \sum_k \int_{\O_k \setminus \O_{k+1}} a^{kp} v dx \\
    &= a^p \sum_{k,j} a^{kp} v(E_j^k),
  \end{split}\end{equation*}

  \noindent since $\O_k \setminus \O_{k+1} = \cup_j E_j^k$ where the sets $E_j^k$ are the sets associated with the family $\{Q_j^k\}$. Then, we obtain

  \begin{equation*}\begin{split}
    \int_{\Rn} \Mm^d(\overrightarrow{f\sigma})^p v dx &\leq a^p \sum_{k,j} \left(\prod_{i=1}^m \frac{1}{|Q_j^k|} \int_{Q_j^k} |f_i|\sigma_i dy_i \right)^p v(E_j^k) \\
    &= a^p \sum_{k,j} v(E_j^k) \left(\prod_{i=1}^m \frac{\sigma_i(Q_j^k)}{|Q_j^k|}\right)^p \left(\prod_{i=1}^m \frac{1}{\sigma(Q_j^k)} \int_{Q_j^k} |f_i|\sigma_i dy_i \right)^p \\
    &= a^p \sum_{Q\in\mathcal{D}} a_Q \left(\prod_{i=1}^m \frac{1}{\sigma(Q_j^k)} \int_{Q_j^k} |f_i| \sigma_i dy_i \right)^p,
  \end{split}\end{equation*}
  where $a_Q = v(E(Q)) \left(\prod_{i=1}^m \frac{\sigma_i(Q)}{|Q|}\right)^p$, if $Q=Q_j^k$ for some $(k,j)$ where $E(Q)$ denotes the corresponding set $E_j^k$ associated to $Q_j^k$, and $a_Q=0$ otherwise.  If we apply the Carleson embedding to these $a_Q$, we will find the desired result provided that

  \begin{equation*}
    \sum_{Q\subset R} a_Q \leq A \int_R \prod_{i=1}^m \sigma_i^{\frac{p}{p_i}}dx, \, R\in\mathcal{D}.
  \end{equation*}

  For $R\in\mathcal{D}$, we obtain

  \begin{equation*}\begin{split}
    \sum_{Q\subset R} a_Q &= \sum_{Q_j^k\subset R} v(E_j^k) \left(\prod_{i=1}^m \frac{\sigma_i(Q_j^k)}{|Q_j^k|}\right)^p \\
    &= \sum_{Q_j^k\subset R} \int_{E_j^k} \left(\prod_{i=1}^m \frac{\sigma_i(Q_j^k)}{|Q_j^k|}\right)^p v(x) dx \\
    &\leq \sum_{Q_j^k\subset R} \int_{E_j^k}  (\Mm(\overrightarrow{\sigma\chi_R}))^p v(x) dx \\
    &\leq [v,\overrightarrow w]_{S_{\overrightarrow P}}^p \prod_{i=1}^m \sigma_i(R)^{\frac{p}{p_i}} \\
    &\leq [v,\overrightarrow w]_{S_{\overrightarrow P}}^p [\overrightarrow \omega]_{RH_{\overrightarrow P}} \int_R \prod_{i=1}^m \sigma_i^{\frac{p}{p_i}}dx,
  \end{split}\end{equation*}
  where in the next to last inequality we have used the $S_{\overrightarrow P}$ condition and in the last inequality we have used the $RH_{\overrightarrow P}$ condition. Thus, by Lemma \ref{Carleson_lemma_multi} we get the desired result and the proof is complete. 

\end{proof}

%
%

\begin{proof}[Theorem \ref{Bp_theorem}]

  To prove this result we proceed using the standard argument as before. We obtain

  \begin{equation*}\begin{split}
    \int_{\Rn} \Mm^d(\overrightarrow{f\sigma})^p v dx &\leq a^p \sum_{k,j} \left(\prod_{i=1}^m \frac{1}{|Q_j^k|} \int_{Q_j^k} |f_i| \sigma_i dy_i \right)^p v(Q_j^k) \\
    &= a^p \sum_{k,j} v(Q_j^k) \left(\prod_{i=1}^m \frac{\sigma_i(Q_j^k)}{|Q_j^k|}\right)^p \left(\prod_{i=1}^m\frac{1}{\sigma_i(Q_j^k)} \int_{Q_j^k} |f_i| \sigma_i dy_i \right)^p \\
    &\leq a^p [v,\overrightarrow{\sigma}]_{B_{\overrightarrow P}} \sum_{k,j} |Q_j^k| \exp\left(\frac{1}{|Q_j^k|}\int_{Q_j^k} \log{\prod_{i=1}^m \sigma_i^{\frac{p}{p_i}}}dx\right) \\ &\times\left(\prod_{i=1}^m\frac{1}{\sigma_i(Q_j^k)} \int_{Q_j^k} |f_i| \sigma_i dy_i \right)^p
  \end{split}\end{equation*}

  And it follows

  \begin{equation*}\begin{split}
    \int_{\Rn} \Mm^d(\overrightarrow{f\sigma})^p v dx &\leq a^p [v,\overrightarrow{\sigma}]_{B_{\overrightarrow P}} \sum_{Q\in\mathcal{D}} a_Q \left(\prod_{i=1}^m\frac{1}{\sigma_i(Q_j^k)} \int_{Q_j^k} |f_i| \sigma_i dy_i \right)^p,
  \end{split}\end{equation*}

  \noindent where in next to last inequality we have used the $B_{\overrightarrow P}$ condition and
  $$a_Q=|Q| \exp\left(\frac{1}{|Q|}\int_{Q} \log{\prod_{i=1}^m \sigma_i^{\frac{p}{p_i}}}dx\right),$$
  if $Q=Q_j^k$ for some $(k,j)$ and $a_Q=0$, otherwise.

  \par Next if we apply Carleson embedding to these $a_Q$, we obtain that \eqref{Bp_eqn} holds provided that

  \begin{equation*}
    \sum_{Q\subset R} a_Q \leq A \int_R \prod_{i=1}^m \sigma_i^{\frac{p}{p_i}}dx, \, R\in\mathcal{D}.
  \end{equation*}

  For $R\in\mathcal{D}$, we have

  \begin{equation*}\begin{split}
    \sum_{Q\subset R} a_Q &\leq \sum_{Q_j^k\subset R} |Q_j^k| \exp\left(\frac{1}{|Q_j^k|}\int_{Q_j^k} \log{\prod_{i=1}^m \sigma_i^{\frac{p}{p_i}}}dx\right) \\
    &\leq 2 \sum_{Q_j^k} |E_j^k|  \exp\left(\frac{1}{|Q_j^k|}\int_{Q_j^k} \log{\prod_{i=1}^m \sigma_i^{\frac{p}{p_i}}}dx\right) \\
    &\leq 2 \sum_{Q_j^k\subset R} \int_{E_j^k} M_0\left(\prod_{i=1}^m \sigma_i^{\frac{p}{p_i}}\chi_R \right)dx \\
    &\leq 2 \int_{\Rn} M_0\left(\prod_{i=1}^m \sigma_i^{\frac{p}{p_i}}\chi_R \right)dx \\
    &\leq 2 e \int_{R} \prod_{i=1}^m \sigma_i^{\frac{p}{p_i}},
  \end{split}\end{equation*}
  where $M_0$ is the (dyadic) logarithmic maximal function described in \cite[Lemma 2.1]{HP} and also discussed in \cite{YM}.  Here we have used that $M_0$ is bounded from $L^1$ into itself, and this concludes the proof of \eqref{Bp_eqn}.
\end{proof}


\begin{proof}[Theorem \ref{Mixed_Ap_Ainfty_theorem}]

Proceeding as we did in the previous theorems, we obtain

\begin{equation*}\begin{split}
  \int_{\Rn} \Mm^d(\overrightarrow{f\sigma})^p v dx &\leq a^p \sum_{k,j} v(Q_j^k) \left(\prod_{i=1}^m \frac{\sigma_i(Q_j^k)}{|Q_j^k|}\right)^p \left(\prod_{i=1}^m\frac{1}{\sigma_i(Q_j^k)} \int_{Q_j^k} |f_i| \sigma_i dy_i \right)^p \\
  &\leq a^p [v,\overrightarrow w]_{A_{\overrightarrow P}} \sum_{Q\in\mathcal{D}} a_Q \left(\prod_{i=1}^m\frac{1}{\sigma_i(Q_j^k)} \int_{Q_j^k} |f_i| \sigma_i dy_i \right)^p,
\end{split}\end{equation*}
where we have used the $A_{\overrightarrow P}$ condition and we have denoted by $a_Q$ the following numbers

\begin{equation*}
  a_Q = \prod_{i=1}^m \sigma_i(Q)^{\frac{p}{p_i}},
\end{equation*}
if $Q=Q_j^k$ for some $(j,k)$, and $a_Q=0$, otherwise.  Therefore it suffices to check that \eqref{Carleson_assumption} holds for every $R\in\mathcal{D}$.  Indeed,

\begin{equation*}\begin{split}
  \sum_{Q\subset R} a_Q &= \sum_{Q_j^k\subset R}\prod_{i=1}^m \sigma_i(Q_j^k)^{\frac{p}{p_i}} = \sum_{Q_j^k\subset R} \prod_{i=1}  \left(\frac{\sigma_i(Q_j^k)}{|Q_j^k|}\right)^{\frac{p}{p_i}}|Q_j^k| \\
  &\leq 2 \sum_{Q_j^k\subset R} |E_j^k| \prod_{i=1}  \left(\frac{\sigma_i(Q_j^k)}{|Q_j^k|}\right)^{\frac{p}{p_i}} \\
  &\leq 2 \sum_{Q_j^k\subset R} \int_{E_j^k} \prod_{i=1}^m M(\sigma_i\chi_R)^{\frac{p}{p_i}} dx \\
  &\leq 2 \int_R \prod_{i=1}^m M(\sigma_i\chi_R)^{\frac{p}{p_i}} dx \\
  &\leq 2 [\overrightarrow{\sigma}]_{W_{\overrightarrow P}^{\infty}} \int_R \prod_{i=1}^m \sigma_i^{\frac{p}{p_i}} dx,
\end{split}\end{equation*}
where $E_j^k$ are the sets associated with the cubes $Q_j^k$ and we have used that $\overrightarrow{\sigma} \in W_{\overrightarrow P}^{\infty}$. Therefore we have proved \eqref{Mixed_eqn}.

\end{proof}

\thanks{\textbf{Acknowledgement.}{This paper was completed while the first author was at the Institute of Mathematics of the University of Seville (I.M.U.S.), Spain. He would like to express his gratitude for the hospitality received there. The authors would like to thank Prof. Carlos P\'erez for his advice and support and his helpful remarks.
\par  \noindent The first author is supported by the National Natural Science Foundation of China (Grant No. 11101353), the Natural Science Foundation of Jiangsu Education Committee (Grant No. 11KJB110018) and the Natural Science Foundation of Jiangsu Province (Grant No. BK2012682). The second author is supported by Junta de Andaluc\'ia (Grant No. P09-FQM-4745).
}

%
%

\end{document}